\documentclass[a4paper,12pt]{amsart}
\usepackage{amssymb,amscd,amsmath,latexsym}
\usepackage[dvips]{graphicx}

\theoremstyle{plain}

\newtheorem{theo}{Theorem}[section]
\newtheorem{prop}[theo]{Proposition}
\newtheorem{lemm}[theo]{Lemma}

\theoremstyle{definition}

\theoremstyle{remark}
\newtheorem*{rema}{Remark}

\numberwithin{equation}{section}

\newcommand{\field}[1]{\mathbb{#1}}
\newcommand{\C}{\field{C}}
\newcommand{\R}{\field{R}}
\newcommand{\Z}{\field{Z}}

\def\T{\mathcal T}

\DeclareMathOperator{\Hom}{Hom}

\DeclareMathOperator{\Ann}{Ann}

\def\S1{\C^*}
\def\T{\mathcal T}
\def\v{\mathcal V}

\begin{document}

\title{Equivariant cohomology distinguishes 
toric manifolds}
\date{\today}

\author[M.~Masuda]{Mikiya Masuda}
\address{Department of Mathematics, Osaka City University, Sugimoto, 
Sumiyoshi-ku, Osaka 558-8585, Japan}
\email{masuda@sci.osaka-cu.ac.jp}

\begin{abstract}
The equivariant cohomology of a space with a group action 
is not only a ring but also an algebra 
over the cohomology ring of the classifying space of the acting group. 
We prove that toric manifolds (i.e. compact smooth toric varieties) are 
isomorphic as varieties if and only if their equivariant cohomology algebras 
are weakly isomorphic. 
We also prove that quasitoric manifolds, which 
can be thought of as a topological counterpart to toric manifolds, 
are equivariantly homeomorphic if and only if their 
equivariant cohomology algebras are isomorphic. 
\end{abstract}

\maketitle

\section{Introduction}

Let $T$ be a $\C^*$-torus of rank $n$, i.e., $T=(\C^*)^n$. 
A toric variety $X$ of complex dimension $n$ is a normal complex algebraic 
variety with an action of $T$ having an open dense orbit.  
A fundamental result in the theory of toric varieties says that there 
is a one-to-one correspondence between toric varieties and fans, and 
among toric varieties, 
compact smooth toric varieties, which we call {\em toric manifolds}, 
are well studied, see \cite{fult93}, \cite{oda88}.  

Suppose two toric manifolds $X$ and $X'$ are isomorphic as varieties.  
Then they are not necessarily equivariantly isomorphic as varieties, but 
\emph{weakly equivariantly} isomorphic as varieties, i.e. 
there is a variety isomorphism $\phi\colon X\to X'$ together with 
an automorphism $\gamma$ of $T$ such that $\phi(tx)=\gamma(t)
\phi(x)$ for any $t\in T$ and $x\in X$. 
This is well-known and follows from the fact that 
the automorphism group of a toric manifold is a linear algebraic group 
with the acting group $T$ as a maximal algebraic torus 
(\cite[Section 3.4]{oda88}).  Therefore, classifying toric manifolds 
up to variety isomorphism is same as that up to weakly equivariant variety 
isomorphism.  

The equivariant cohomology of a toric variety $X$ is by definition 
\[
H^*_T(X):=H^*(ET\times_T X)
\]
where $ET$ is the total space of the universal principal $T$-bundle and 
$ET\times_T X$ is the oribit space of $ET\times X$ by the diagonal $T$-action.
$H^*_T(X)$ contains a lot of geometrical 
information on $X$, but its ring structure does not reflect enough 
geometrical information on $X$. 
In fact, when $X$ is a toric manifold, $H^*_T(X)$ as a ring 
is the face ring of the underlying simplicial complex $\Sigma$ 
of the fan of $X$ and determined by $\Sigma$. 
There are toric manifolds which are not isomorphic as varieties but 
have the same underlying simplicial complex, so equivariant 
cohomology as a ring does not distinguish toric manifolds.  

However, $H^*_T(X)$ is not only a ring but also 
an algebra over $H^*(BT)$ through the projection map from 
$ET\times_T X$ onto $ET/T=BT$. This algebra structure contains more 
geometrical information on $X$.  
If two toric manifolds $X$ and $X'$ are isomorphic as varieties, then they 
are weakly equivariantly isomorphic as varieties as remarked above, so that 
$H^*_T(X)$ and $H^*_T(X')$ are \emph{weakly isomorphic} as 
algebras over $H^*(BT)$, i.e., there is a ring isomorphism $\Phi\colon 
H^*_T(X')\to H^*_T(X)$ together with an automorphism $\gamma$ of $T$ such that 
$\Phi(u\omega)=\gamma^*(u)\Phi(\omega)$ for any $u\in H^*(BT)$ and $\omega
\in H^*_T(X')$ where $\gamma^*$ denotes the automorphism of $H^*(BT)$ induced 
by $\gamma$. Our main result asserts that the converse holds. 

\begin{theo} \label{main}
Two toric manifolds are (weakly equivariantly) isomorphic as varieties if 
and only if their equivariant cohomology algebras 
are weakly isomorphic. 
\end{theo}

The theorem above leads us to ask how much information ordinary 
cohomology contains for toric manifolds, in particular we may ask  
whether two toric manifolds are homeomorphic (or diffeomorphic) 
if their ordinary cohomology rings are isomorphic. 
The question is affirmatively solved in some cases 
(\cite{ma-pa06}, \cite{ch-ma-su07}) and the author does not know any 
counterexample.  

This paper is organized as follows.  In Section 2 we review how the 
equivariant cohomology of a toric manifold $X$ is related to the fan of $X$, 
and prove Theorem~\ref{main} in Section 3.  In Section 4 we observe 
that our argument also works with some modification for quasitoric 
manifolds which are a topological counterpart to toric manifolds. 

\medskip
\noindent
{\bf Acknowledgment.}  
The first version of this paper was written while the author was visiting 
Fudan University in fall 
2006.  He would like to thank Fudan University and Zhi L\"u 
for the invitation and providing excellent circumstances to do research. 
He also would like to thank Taras Panov, Hiroshi Sato and Dong Youp Suh 
for useful discussions. 

\bigskip

\section{Equivariant cohomology and fan} \label{eqfan}

Throughout this and next sections, $X$ will denote a toric manifold of 
complex dimension $n$ unless otherwise stated. 
In this section we shall review how the equivariant cohomology 
of $X$ is related to the fan of $X$. 
The reader will find that most of the arguments in this and next sections 
work with a compact torus $(S^1)^n$ instead of $T=(\C^*)^n$. 

There are only finitely many $T$-invariant divisors in $X$, which we denote 
by $X_1,\dots,X_m$.  Each $X_i$ is a complex codimension-one invariant 
closed submanifold of $X$ and fixed pointwise by some $\C^*$-subgroup 
of $T$.  Since $X$ and $X_i$ are complex manifolds, 
they have canonical orientations.  Let $\tau_i\in H^2_T(X)$ be the 
Poincar\'e dual of $X_i$ viewed as an equivariant cycle in $X$, 
in other words, $\tau_i$ is the image of the unit $1\in H^0_T(X_i)$ 
by the equivariant Gysin homomorphism $\colon H^0_T(X_i)\to H^2_T(X)$ induced 
from the inclusion map $\colon X_i\to X$. 
We call $\tau_i$ the {\em Thom class} of $X_i$. 

We abbreviate a set $\{1,\dots,m\}$ as $[m]$. 
The invariant divisors $X_i$ intersect transversally, so a cup product 
$\prod_{i\in I}\tau_i$ for a subset $I$ of $[m]$ is the Poincar\'e 
dual of the intersection $\cap_{i\in I}X_i$. 
In particular, $\prod_{i\in I}\tau_i=0$ if $\cap_{i\in I}X_i=\emptyset$. 
Since $H^*(X)$ is generated by elements in $H^2(X)$ as a ring 
(see \cite[section 3.3]{oda88}), 
we see that $H^*_T(X)$ is generated by $\tau_i$'s as a ring and there is 
no more relation among $\tau_i$'s than those mentioned above, 
see {\cite[Proposition 3.4]{masu99}} for example.  Namely we have 

\begin{prop}\label{Tring}
\[
H^*_T(X)=\Z[\tau_1,\dots,\tau_m]/(\prod_{i\in I}\tau_i \mid \bigcap_{i\in I}X_i
=\emptyset) \quad \text{as ring}
\]
where $I$ runs all subsets of $[m]$ such that $\bigcap_{i\in I}X_i
=\emptyset$. 
\end{prop}

We set  
\[
\Sigma:=\{ I\subset [m]\mid \bigcap_{i\in I}X_i\not=\emptyset\}.
\]
This is an abstract simplicial complex of dimension $n-1$ and the proposition 
above says that $H^*_T(X)$ is the face ring (or Stanley-Reisner ring) of the 
simplicial complex $\Sigma$. 

Let $\pi\colon ET\times_T X\to ET/T=BT$ be the projection on the first factor. 
Through $\pi^*\colon H^*(BT)\to H^*_T(X)$, one can regard $H^*_T(X)$ as 
an algebra over $H^*(BT)$.  Since $T$ is a torus of rank $n$, 
$H^*(BT)$ is a polynomial ring in $n$ variables of degree two, in particular, 
it is generated by elements of degree two as a ring.  Therefore, one can 
find the algebra structure of $H^*_T(X)$ over $H^*(BT)$ if one knows 
how elements in $H^2(BT)$ map to $H^2_T(X)$ by $\pi^*$. 

\begin{prop}
\label{Talge}
To each $i\in [m]$, there is a unique element $v_i\in H_2(BT)$ 
such that 
\begin{equation} \label{keyeq}
\pi^*(u)=\sum_{i=1}^m\langle u,v_i\rangle \tau_i\quad \text{for any $u\in 
H^2(BT)$}
\end{equation}
where $\langle\ ,\ \rangle$ is the natural pairing between cohomology and 
homology. 
\end{prop}

\begin{rema}  The identity (\ref{keyeq}) corresponds to the identity 
in \cite[Lemma in p.61]{fult93} in algebraic geometry, 
which describes a principal divisor 
as a linear combination of the $T$-invariant divisors $X_i$. 
\end{rema}

\begin{proof} The proposition is proved in 
\cite[Lemma 9.3]{ha-ma03} and \cite[Lemma 1.5]{masu99}.  
But for the reader's convenience we shall reproduce the proof given 
in \cite[Lemma 9.3]{ha-ma03}.  
By Proposition~\ref{Tring} $H^2_T(X)$ is freely generated by $\tau_1, 
\dots,\tau_m$ over $\Z$.  Therefore, for each $u\in H^2(BT)$, one can 
uniquely express $\pi^*(u)\in H^2_T(X)$ as 
\[
\pi^*(u)=\sum_{i=1}^mv_i(u)\tau_i
\]
with integers $v_i(u)$ depending on $u$.  We view $v_i(u)$ as a function 
of $u$. Since $\pi^*$ is a homomorphism, 
the function $v_i(u)$ is linear; so there 
is a unique $v_i\in H_2(BT)$ such that $v_i(u)=\langle u,v_i\rangle$. 
\end{proof}

The vectors $v_i$ have a nice geometrical meaning, which we shall explain. 
The group $\Hom(\S1,T)$ of homomorphisms 
from $\S1$ to $T$ can be identified with $H_2(BT)$ as follows. 
An element $\rho$ of $\Hom(\S1,T)$ induces a continuous map 
$\bar\rho\colon B\S1\to BT$ between classifying spaces 
and $H_2(B\S1)$ is isomorphic to $\Z$; so 
once we choose and fix a generator, say $\alpha$, of $H_2(B\S1)$, we get 
an element $\bar\rho_*(\alpha)\in H_2(BT)$.  A correspondence $\colon 
\rho \to \bar\rho_*(\alpha)$ gives an isomorphism from 
$\Hom(\S1,T)$ to $H_2(BT)$ and we denote by $\lambda_v$ the element of 
$\Hom(\S1,T)$ corresponding to $v\in H_2(BT)$. 
It turns out that $\lambda_{v_i}(\S1)$ 
is the $\C^*$-subgroup of $T$ fixing $X_i$ pointwise, 
see \cite[Lemma 1.10]{masu99} for example. 

We have obtained two data from $X$, one is the abstract simplicial complex 
$\Sigma$ and the other is the set of vectors $v_1,\dots,v_m$ in $H_2(BT)$. 
To each $I\in \Sigma$ we form a cone in $H_2(BT)\otimes\R=H_2(BT;\R)$ 
spanned by $v_i$'s $(i\in I)$. Then the collection of these cones is the fan 
of $X$. Precisely speaking, we need to add the $0$-dimensional cone consisting 
of the origin to this collection to satisfy 
the conditions required in the definition of fan, see \cite{fult93} or 
\cite{oda88}. 
The $0$-dimensional cone corresponds to the empty subset of $[m]$. 
Although we formed cones using the data $\Sigma$ and $\{v_i\}$ to define 
the fan of $X$, we may think of a pair $(\Sigma,\{v_i\})$ as the fan of $X$. 

As is well known $X$ can be recovered 
from the fan of $X$. There are at least three ways (gluing affine spaces, 
taking quotient by a $\C^*$-torus or symplectic reduction) to recover $X$ 
from the fan of $X$. 
We shall recall the quotient construction. 
For $x=(x_1,\dots,x_m)\in \C^m$ we define $I(x)=\{ i\mid x_i=0\}$. 
We note that $(\C^*)^m$ acts on $\C^m$ via coordinatewise scalar 
multiplication. 

\begin{prop}[see \cite{cox95}] \label{recov}
Let $X$ be a toric manifold and $(\Sigma,\{v_i\})$ be the fan of $X$. 
We consider 
\[
Y:=\{ x\in \C^m\mid I(x)\in \Sigma\cup\{\emptyset\}\}
\]
and a homomorphism 
\[
\v\colon (\C^*)^m\to (\C^*)^n=T
\]
defined by 
$$\v(g_1,\dots,g_m)=\prod_{i=1}^m \lambda_{v_i}(g_i).$$ 
Then $Y$ is invariant under the $(\C^*)^m$-action, 
the kernel $\ker \v$ of $\v$ acts on $Y$ freely and 
the quotient $Y/\ker \v$ with the induced $T$-action 
is a toric manifold equivariantly isomorphic to $X$. 
\end{prop} 

\bigskip

\section{Poof of Theorem~\ref{main}}
We continue to use the notation in Section 2. 
Let $X^T$ denote the set of $T$-fixed points in $X$.  As is well known, 
it consists of finitely many points. 
For $\xi\in H^2_T(X)$, we denote its restriction to $p\in X^T$ by $\xi|p$ 
and define 
\[
Z(\xi):=\{ p\in X^T\mid \xi|p=0\}.
\]

\begin{lemm} \label{length}
Express $\xi=\sum_{i=1}^m a_i\tau_i$ with integers 
$a_i$. If $a_i\not=0$ for some $i$, then $Z(\xi)\subset Z(\tau_i)$.  Moreover, 
if $a_i\not=0$ and $a_j\not=0$ for some different $i$ and $j$, 
then $Z(\xi)\subsetneq Z(\tau_i)$.
\end{lemm}

\begin{proof}
Let $p\in X^T$. 
Since $\tau_i$ is the Poincar\'e dual of $X_i$ viewed as an equivariant cycle 
in $X$, $\tau_i|p=0$ if $p\notin X_i$.  
Moreover, if $p\in X_i$, then $\tau_i|p\in H^2_T(p)=H^2(BT)$ is the 
equivariant Euler class of 
the complex one-dimensional normal $T$-representation at $p$ to $X_i$. 
This implies that 
\begin{equation} \label{taup}
\tau_i|p=0 \quad\text{if and only if}\quad p\notin X_i
\end{equation} 
and that there are exactly $n$ number of $X_i$'s 
containing $p$ and $\{\tau_i|p\mid p\in X_i\}$ forms a basis of $H^2(BT)$. 

Suppose $p\in Z(\xi)$. Then $0=\xi|p=\sum_{i=1}^m a_i\tau_i|p$ and 
it follows from the observation above that 
$\tau_i|p=0$ if $a_i\not=0$.  This proves the former statement in the lemma. 

If both $a_i$ and $a_j$ are non-zero, then $Z(\xi)\subset Z(\tau_i)\cap 
Z(\tau_j)$ by the former statement in the lemma.  
Therefore, it suffices to prove 
that $Z(\tau_i)\cap Z(\tau_j)$ is properly contained in $Z(\tau_i)$. 
Suppose that 
$Z(\tau_i)\cap Z(\tau_j)=Z(\tau_i)$.  Then $Z(\tau_j)\supset Z(\tau_i)$, 
so $X_j^T\subset X_i^T$ by (\ref{taup}). 
This implies that $X_j=X_i$, a contradiction. 
\end{proof}

Let $S=H^*(BT)\backslash \{0\}$ and let $S^{-1}H^*_T(X)$ denote the localized 
ring of $H^*_T(X)$ by $S$. 
Since $H^{odd}(X)=0$, $H^*_T(X)$ is free as a module over $H^*(BT)$. 
Hence the natural map 
\[
H^*_T(X)\to S^{-1}H^*_T(X)\cong S^{-1}H^*_T(X^T)=
\bigoplus_{p\in X^T}S^{-1}H^*_T(p)
\]
is injective, where the isomorphism above is induced from the 
inclusion map from $X^T$ to $X$ and is a consequence of the Localization 
Theorem in equivariant cohomology (\cite[p.40]{hsia75}).  The annihilator 
$$\Ann(\xi):=\{\eta \in S^{-1}H^*_T(X)\mid \eta\xi=0\}$$ 
of $\xi$ in $S^{-1}H^*_T(X)$ is nothing but sum of $S^{-1}H^*_T(p)$ 
over $p$ with $\xi|p=0$.  Therefore it is a free $S^{-1}H^*(BT)$ 
module of rank $|Z(\xi)|$.  
Since $\Ann(\xi)$ is defined using the algebra structure of $H^*_T(X)$, 
$|Z(\xi)|$ is an invariant of $\xi$ 
depending only on the algebra structure of $H^*_T(X)$.  
We note that $|Z(\xi)|$ is invariant under an algebra isomorphism. 
We call $|Z(\xi)|$ the {\em zero-length} of $\xi$. 

\begin{lemm} \label{Thom}
Let $X'$ be another toric manifold ($X'$ might be same as $X$). 
If $f\colon H^*_T(X)\to H^*_T(X')$ is an algebra isomorphism, 
then $f$ maps the Thom classes in $H^2_T(X)$ to the Thom classes in 
$H^2_T(X')$ bijectively up to sign.  
\end{lemm}

\begin{proof}
We classify the Thom classes according to their zero-length. 
Let $\T_1$ be the set of Thom classes in $H^2_T(X)$ with largest 
zero-length, and let $\T_2$ be the set of Thom classes in $H^2_T(X)$ with 
second largest zero-length, and so on.  Similarly we define $\T_1', \T_2'$ 
and so on for the Thom classes in $H^2_T(X')$.  

Let $m_k$ (resp. $m_k'$) be the zero-length of elements in $\T_k$ (resp. 
$\T_k'$).  
Since both $f$ and $f^{-1}$ preserve zero-length and are isomorphisms, 
$m_1=m_1'$ and $f$ maps $\T_1$ to $\T_1'$ bijectively up to sign by 
Lemma~\ref{length}.  
Take an element $\tau_i$ from $\T_2$. Since $\T_1$ and 
$\T_1'$ are preserved under $f$ and $f^{-1}$, $f(\tau_{i})$ is 
not a linear combination of elements in $\T_1'$ .  
This together with Lemma~\ref{length} means 
that $m_2\leq m_2'$.  The same argument for $f^{-1}$ instead of $f$ shows 
that $m_2'\leq m_2$, so that $m_2=m_2'$.  Again, this together with 
Lemma~\ref{length} implies that $f$ maps $\T_2$ to $\T_2'$ bijectively 
up to sign.  The lemma follows by repeating this argument.  
\end{proof}

Now we shall complete the proof of Theorem~\ref{main}. 
Let $X$ and $X'$ be two toric manifolds whose equivariant cohomology algebras 
over $H^*(BT)$ are weakly isomorphic. We note that changing the action of $T$ 
on $X$ through an automorphism of $T$, we may assume that $H^*_T(X)$ and 
$H^*_T(X')$ are isomorphic as algebras over $H^*(BT)$. 

We put a prime for notation for $X'$ corresponding to 
the Thom classes $\tau_i$, the abstract simplicial complex $\Sigma$ and 
the vectors $v_i$ etc. for $X$. 
Let $f\colon H^*_T(X)\to H^*_T(X')$ be an isomorphism of algebras over 
$H^*(BT)$. 
By Lemma~\ref{Thom}, the number of the Thom classes in $H^2_T(X)$ is same as 
that in $H^2_T(X')$ and there is a permutation $\bar f$ on $[m]$ such that 
\begin{equation} \label{ftau}
f(\tau_i)=\epsilon_i\tau'_{\bar f(i)} \quad\text{with $\epsilon_i=\pm 1$.}
\end{equation} 
If $I\subset [m]$ is an element of $\Sigma$, then 
$\prod_{i\in I}\tau_i$ is non-zero by Proposition~\ref{Tring} and hence so is 
$f(\prod_{i\in I}\tau_i)=\prod_{i\in I}\epsilon_i\tau'_{\bar f(i)}$.  
Therefore a subset $\{ \bar f(i)\mid i\in I\}$ of $[m]$ is 
a simplex in $\Sigma'$ whenever $I$ is a simplex in $\Sigma$, which means 
that $\bar f$ induces a simplicial map from $\Sigma$ to $\Sigma'$. 
Applying the same argument to the inverse of $f$, we see that the induced 
simplicial map has an inverse, so that it is an isomorphism. 

Since $f$ is an algebra map over $H^*(BT)$, 
${\pi'}^*=f\circ\pi^*$.  Therefore, 
sending the identity (\ref{keyeq}) by $f$ and using (\ref{ftau}), we have 
\[
{\pi'}^*(u)=f(\pi^*(u))=\sum_{i=1}^m \langle u,v_i\rangle f(\tau_i)
=\sum_{i=1}^m \langle u,v_i\rangle \epsilon_i\tau'_{\bar f(i)}.
\]
Comparing this with the identity (\ref{keyeq}) for $X'$ and 
noting that $\bar f$ is a permutation on $[m]$, we conclude that 
\begin{equation} \label{fv}
\epsilon_iv_i=v'_{\bar f(i)}\quad \text{for each $i$.}
\end{equation} 

We identify $\Sigma$ with $\Sigma'$ through the isomorphism induced by 
$\bar f$, 
so that we may think of $\bar f$ as the identity map and then the identity 
(\ref{fv}) turns into 
\[
\epsilon_i v_i=v_i'.
\]
By Proposition~\ref{recov} we may assume $X=Y/\ker\v$ and $X'=Y'/\ker\v'$.  
Since $\Sigma'$ is identified with $\Sigma$, we have $Y=Y'$.  
Therefore it suffices to check that $\ker \v=\ker \v'$. 
Since $\lambda_{-v}(g)=\lambda_{v}(g)^{-1}=\lambda_{v}(g^{-1})$ for 
$v\in H_2(BT)$ and $g\in 
\S1$, an automorphism $\rho$ of $(\C^*)^m$ defined by 
$$\rho(g_1,\dots,g_m)=(g_1^{\epsilon_1},\dots,g_m^{\epsilon_m})$$
satisfies $\v\circ \rho=\v'$. This implies $\ker \v=\ker\v'$ and 
completes the proof of Theorem~\ref{main}.

\bigskip

\section{Quasitoric manifolds}

Davis-Januszkiewicz 
\cite{da-ja91} introduced the notion of what is now called a 
\emph{quasitoric manifold}, see \cite{bu-pa02}.  
A quasi-toric manifold is a closed smooth 
manifold of even dimension, say $2n$, with a smooth action of a compact torus 
group $(S^1)^n$ of dimension $n$ such that the action is locally isomorphic 
to a faithful $(S^1)^n$-representation of real dimension $2n$ and that 
the orbit 
space is combinatorially a simple convex polytope.  A toric manifold with the 
action restricted to the maximal compact toral subgroup of $T$ 
often provides an example of a quasitoric manifold, e.g. this is the 
case when $X$ is projective.  However, there are many quasitoric manifolds 
which do not arise from a toric manifold.  For instance, $\C P^2\#\C P^2$ 
with an appropriate action of $(S^1)^2$ is a quasitoric manifold but 
does not arise from a toric manifold because $\C P^2\#\C P^2$ 
does not allow a complex (even almost complex) structure. 
We note that the equivariant cohomology of a quasitoric manifold of dimension 
$2n$ is an algebra over $H^*(B(S^1)^n$) similarly to the toric case. 
The purpose of this section is to prove the following. 

\begin{theo} \label{quasi}
Two quasitoric manifolds are equivariantly homeomorphic if 
their equivariant cohomology algebras are isomorphic.  
\end{theo}


\begin{proof} 
When $X$ is a quasitoric manifold, we take $X_i$ to be a connected real 
codimension-two closed submanifold of $X$ fixed pointwise by some circle 
subgroup of $(S^1)^n$. Then 
the proof for Theorem~\ref{main} almost works if we replace $\C^*$ by $S^1$ 
(and hence $T=(\C^*)^n$ by $(S^1)^n$). The only problem is that we do 
not have Proposition~\ref{recov} for quasitoric manifolds, so that the 
last paragraph in the previous section needs to be modified. 
In the sequel, it suffices to prove that the existence of an isomorphism 
$\bar f\colon \Sigma\to \Sigma'$ satisfying (\ref{fv}) implies that the two 
quasitoric manifolds $X$ and $X'$ are equivariantly homeomorphic. 

Let $P$ be the orbit space of $X$ by the action of 
$(S^1)^n$ and let $q\colon X\to P$ be the quotient map.  The orbit space 
$P$ is a simple convex polytope by the definition of quasitoric manifold.  
Then $P_i:=q(X_i)$ is a facet (i.e., a codimension-one face) of $P$. 
The dual polytope $P^*$ of $P$ is a simplicial polytope 
and its boundary complex agrees with $\Sigma$. The vertices of $\Sigma$ 
bijectively correspond to the facets of $P$ so that $v_i$ is assigned 
to $P_i$. The vectors $v_i$ form a characteristic function on $P$ 
introduced in \cite{da-ja91}.  Any (proper) face of $P$ is obtained as an 
intersection $P_I:=\cap_{i\in I}P_i$ for some $I\in\Sigma$. 
We define $P_\emptyset$ to be $P$ itself. 
For $I\in \Sigma$ we denote by $S_I$ a subgroup of $(S^1)^n$ 
generated by circle subgroups $\lambda_{v_i}(S^1)$ for $i\in I$. 
We define $S_\emptyset$ to be the unit group. 
Associated with a pair $(P,\{v_i\})$ we form a quotient space 
\[
X(P,\{v_i\}):=P\times (S^1)^n/\sim. 
\]
Here $(p_1,g_1)\sim (p_2,g_2)$ if and only if $p_1=p_2$ and $g_1^{-1}g_2
\in S_I$ where $I\in \Sigma\cup\{\emptyset\}$ is determined by the condition 
that $p_1=p_2$ is contained in the interior of $P_I$.  The natural action 
of $(S^1)^n$ on the product $P\times (S^1)^n$ descends to an action on 
$X(P,\{v_i\})$ and $X$ is equivariantly homeomorphic to 
$X(P,\{v_i\})$ (see \cite[Proposition 1.8]{da-ja91}). 

As before, we put a prime to denote elements for $X'$ corresponding to 
$P, v_i$ and $\Sigma$. 
The isomorphism $\bar f\colon \Sigma\to \Sigma'$ induces an isomorphism 
from $P^*$ to ${P'}^*$ and then a face-preserving homeomorphism 
from $P$ to $P'$ which we 
denote by $\varphi$.  A map $\varphi\times id \colon P\times (S^1)^n\to 
P'\times (S^1)^n$ descends to a map from $X(P,\{v_i\})$ to $X(P',\{v_i'\})$ 
by virtue of (\ref{fv}) and the resulting map is an 
equivariant homeomorphism, so $X$ is equivariantly homeomorphic to $X'$. 
\end{proof}

Similarly to the toric case, 
it would be interesting to ask whether two quasitoric manifolds are 
homeomorphic (or diffeomorphic) if their ordinary cohomology rings are 
isomorphic, see \cite{ma-pa06} and \cite{ch-ma-su07} for some 
partial affirmative solutions. 

\begin{rema} 
Davis-Januszkiewicz \cite{da-ja91} also introduced the notion of a real 
version of quasitoric manifold, which they call a {\em small cover}.  
A small cover is a closed smooth manifold of dimension, say $n$, with a smooth 
action of a rank $n$ mod two torus group $(\Z_2)^n$ such that 
the action is locally isomorphic to a faithful $(\Z_2)^n$-representation of 
real dimension $n$ and that the orbit space is combinatorially a simple 
convex polytope.  Our argument also works for small covers with $\Z_2$ 
coefficient, so that 
small covers are equivariantly homeomorphic if their equivariant cohomology 
algebras with $\Z_2$ coefficient are isomorphic. 
\end{rema}

\end{document}